\documentclass
[superscriptaddress,secnumarabic,amssymb,amsmath,nobibnotes,aps,prd,showkeys,showpacs,nofootinbib,onecolumn]{revtex4}%
\usepackage{graphicx}
\usepackage{epsf}
\usepackage{bm}
\usepackage{amsmath}
\usepackage{amsfonts}
\usepackage{amssymb}%
\providecommand{\U}[1]{\protect\rule{.1in}{.1in}}

\newcommand{\be}{\begin{equation}}
\newcommand{\ee}{\end{equation}}

\newcommand{\mincir}{\raise
-3.truept\hbox{\rlap{\hbox{$\sim$}}\raise4.truept\hbox{$<$}\ }}
\newcommand{\magcir}{\raise
-3.truept\hbox{\rlap{\hbox{$\sim$}}\raise4.truept\hbox{$>$}\ }}

\ifx\pdfoutput\relax\let\pdfoutput=\undefined\fi
\newcount\msipdfoutput
\ifx\pdfoutput\undefined\else
\ifcase\pdfoutput\else
\ifx\paperwidth\undefined\else
\ifdim\paperheight=0pt\relax\else\pdfpageheight\paperheight\fi
\ifdim\paperwidth=0pt\relax\else\pdfpagewidth\paperwidth\fi
\fi\fi\fi
\begin{document}
\title{Nonlocal Representation of the $sl(2,R)$ Algebra for the Chazy equation}
\author{Andronikos Paliathanasis}
\email{anpaliat@phys.uoa.gr}
\affiliation{Instituto de Ciencias F\'{\i}sicas y Matem\'{a}ticas, Universidad Austral de
Chile, Valdivia, Chile}
\affiliation{Institute of Systems Science, Durban University of Technology, PO Box 1334,
Durban 4000, Republic of South Africa}
\author{Sameerah Jamal}
\email{sameerah.jamal@wits.ac.za}
\affiliation{School of Mathematics and Centre for Differential Equations,Continuum
Mechanics and Applications, University of the Witwatersrand, Johannesburg,
South Africa}
\author{P.G.L. Leach}
\email{leach.peter@ucy.ac.cy}
\affiliation{Department of Mathematics and Institute of Systems Science, Research and
Postgraduate Support, Durban University of Technology, POB 1334 Durban 4000,
Republic of South Africa}
\affiliation{School of Mathematics, Statistics and Computer Science, University of
KwaZulu-Natal, Private Bag X54001, Durban 4000, Republic of South Africa}

\begin{abstract}
A demonstration of how the point symmetries of the Chazy Equation become
nonlocal symmetries for the reduced equation is discussed. Moreover we
construct an equivalent third-order differential equation which is related to
the Chazy Equation under a generalized transformation, and find the point
symmetries of the Chazy Equation are generalized symmetries for the new
equation. With the use of singularity analysis and a simple coordinate
transformation we construct a solution for the Chazy Equation which is given
by a Right Painlev\'{e} Series. The singularity analysis is applied to the new
third-order equation and we find that it admits two solutions, one given by a
Left Painlev\'{e} Series and one given by a Right Painlev\'{e} Series where
the leading-order behaviors and the resonances are explicitly those of the
Chazy Equation.

\end{abstract}
\keywords{Chazy Equation; Lie symmetries; Nonlocal symmetries; Singularity analysis}\maketitle
\date{\today}

\section{Introduction}

In 1911 Jean Chazy presented a paper \cite{Chazy11a} in Comptes Rendus in
which he extended the classification of ordinary differential equations
according to their possession of moveable singularities initiated by
Painlev\'{e} \cite{Painleve95a,Painleve00a,Painleve02a,Painleve06a} and his
colleagues \cite{Gambier89a,Garnier12a} to equations of higher order. Of the
third-order equations in Chazy's catalogue, the one labelled 12a, namely
\begin{equation}
y^{\prime\prime\prime}-2yy^{\prime\prime}+3y^{\prime2}=0, \label{eq.01}%
\end{equation}
has the peculiar property that it possesses a `natural barrier' through which
the solution does not pass.

In the spirit of singularity analysis initiated by Sophie Kowalevskaya
\cite{Kowalevski88a} and enshrined in the ARS Algorithm
\cite{Ablowitz78a,Ablowitz80a, Ablowitz80b} the moveable singularity is a
simple pole and the resonances are $-1$, $-2$ and $-3$ (see below for
details). Conventional wisdom had it that (\ref{eq.01}) did not possess the
Painlev\'{e} Property although considerable ingenuity was applied to resolve
the conundrum of negative resonances \cite{Fordy91a,Conte91a}. A hint of the
explanation is to be found in \cite{Lemmer93a} and in \cite{Feix97a}, but the
full explanation had to wait until 2006 when Andriopoulos \textit{et al}
\cite{Andriopoulos06a} showed that the expansion about the singularity did not
have to be in the vicinity of the singularity, but could be over an annulus
centred on the singularity or the exterior of a disc also centred on the
singularity. The first is found when (apart from the generic $-1$) the
resonances are positive, the second when the resonances are of mixed sign and
the third when they are negative. This becomes very obvious when one thinks in
terms of the three possible forms which a Laurent Expansion can take.

We are concerned with the properties of (\ref{eq.01}) and equations derived
from it in terms of both symmetry and singularity analyses. We recall that the
solution of (\ref{eq.01}) was given in terms of hypergeometric functions by
Olver and Clarkson \cite{Olver98a} as the solutions of a second-order
equation, see also \cite{abl1}. Specifically, in \cite{abl1} it was proven
that the Chazy equation can be derived\ analytically from Schwarz triangle
functions for specific values of the free parameters. Moreover, a new solution
was derived by using the Ramanujan approach for the Eisenstein series.

In this work we look at a different route to the reduction of the Chazy
equation. In particular, we reduce equation (\ref{eq.01}) to a nonautonomous
second-order equation which does not possess any Lie Point Symmetries, then
render the equation in autonomous form by increasing the order to the third by
means of a different transformation than that originally used to reduce
(\ref{eq.01}), and lastly investigate its singularity properties.

\section{Symmetries and Lie invariants}

For the Chazy Equation, (\ref{eq.01}), the Lie point symmetry vectors are
given to be \cite{Olver98a}
\begin{equation}
X_{1}=\partial_{x}~~,~~X_{2}=x\partial_{x}-y\partial_{y}%
\end{equation}%
\begin{equation}
~X_{3}=x^{2}\partial_{x}-\left(  2xy+6\right)  \partial_{y}%
\end{equation}
and constitute a representation of the $sl\left(  2,R\right)  $ algebra.

Reduction with $X_{1}$ gives the differential invariants
\begin{equation}
w=y^{\prime}~,~r=y. \label{eq.00a}%
\end{equation}

The reduced second-order differential equation,%
\begin{equation}
\frac{d^{2}w}{dr^{2}}+\frac{1}{w}\left(  \frac{dw}{dr}\right)  ^{2}-2\frac
{r}{w}\left(  \frac{dw}{dr}\right)  +3=0, \label{eq.02}%
\end{equation}
admits the point symmetry vector,$~r\partial_{r}+2w\partial_{w},~$which is the
original symmetry vector $X_{2}.~$However, under the same transformation
(\ref{eq.00a}), the symmetry vector $X_{3}$ becomes%
\begin{equation}
X_{3}=\left(  6+2r\int\frac{dr}{w}-w\left(  \int\frac{dr}{w}\right)
^{2}\right)  \partial_{r}+\left(  2r+4w\int\frac{dr}{w}\right)  \partial_{w}
\label{eq.03}%
\end{equation}
and is a nonlocal symmetry for the reduced equation (\ref{eq.02}).%
\[
w=\frac{dy}{dx}~,~r=y
\]%
\[
X_{3}=x^{2}\partial_{x}-\left(  2xy+6\right)  \partial_{y}%
\]

From the (reduced) symmetry $X_{2}$ we find the invariants%
\begin{equation}
s=\frac{w}{r^{2}}~,~\phi=\frac{1}{r}\frac{dw}{dr}. \label{eq.03a}%
\end{equation}
Their application to (\ref{eq.02}) gives the first-order differential equation
of Abel type,
\begin{equation}
s\left(  \phi-2s\right)  \frac{d\phi}{ds}+\phi^{2}+\left(  s-2\right)
\phi+3s=0. \label{eq.04}%
\end{equation}

Moreover, when we use (\ref{eq.03a}), the nonlocal symmetry$~$(\ref{eq.03})
becomes an exponential nonlocal symmetry for (\ref{eq.04}) and can be used to
reduce the latter equation to an algebraic equation which gives that the
solution is expressed in terms of the hypergeometric functions \cite{Olver98a}.

Up to this point we have applied the reduction with the symmetry vector
$X_{1}$. It is well known that reduction with $X_{3}$ provides the same
results as $X_{1}$. Therefore we continue with the application of the Lie
invariants of $X_{2}$ which are
\begin{equation}
\chi=xy~,~\psi=x^{2}y^{\prime}. \label{eq.05}%
\end{equation}
We consider $\psi$ to be the new dependent variable and $\chi$ the new
independent variable. Hence the reduced equation is the nonautonomous
second-order differential equation,%
\begin{equation}
\left(  \chi+\psi\right)  ^{2}\frac{d^{2}\psi}{d\chi^{2}}+\left(  \chi
+\psi\right)  \left(  \frac{d\psi}{d\chi}-2\left(  2+\chi\right)  \right)
\frac{d\psi}{d\chi}+\left(  6+4\chi+3\psi\right)  \psi=0. \label{eq.06}%
\end{equation}

It is straightforward to see that equation (\ref{eq.06}) does not admit any
Lie point symmetries. However, in this case the symmetries, $X_{1}$ and
$X_{3}$, become exponential nonlocal symmetries and are%
\begin{equation}
\bar{X}_{1}=e^{-\int\left(  \chi+\psi\right)  ^{-1}d\chi}\left(  \chi
+\psi\right)  \partial_{\chi}~~
\end{equation}
and
\begin{equation}
\bar{X}_{3}=e^{\int\left(  \chi+\psi\right)  ^{-1}d\chi}\left[  \left(
\psi-\chi-6\right)  \partial_{\chi}-2\left(  \chi+\psi\right)  \partial_{\psi
}\right]  .
\end{equation}

The symmetry $X_{2}$ has not disappeared, but it has changed form and does not
play any role in the reduced equation. As far as the second-order differential
equation is concerned, (\ref{eq.06}), we apply the transformation%
\[
\chi=\Phi\left(  s\right)  ~,~\psi=\frac{d\Phi\left(  s\right)  }{ds}%
\]
whereby the scond-order equation can be written as the autonomous third-order
equation
\begin{align}
0  &  =\Phi_{,s}\left(  \Phi_{,s}+\Phi\right)  ^{2}\Phi_{,sss}+\Phi
_{,s}\left(  3\Phi_{,s}+6+4\Phi\right) \nonumber\\
&  -\Phi_{,ss}\left(  \Phi_{,s}+\Phi\right)  \left(  \Phi\left(  \Phi
_{,ss}+2\Phi_{,s}^{2}\right)  +4\Phi_{,s}^{4}\right)  \label{eq.07}%
\end{align}
in which the only admitted Lie point symmetry is the autonomous symmetry,
$\partial_{s}$. However, now the nonlocal symmetries, $\bar{X}_{1}~$and
$\bar{X}_{3}$, are changed in form and are generalized symmetries, that is,%
\begin{equation}
\bar{X}_{1}=\left(  \frac{\Phi_{,ss}}{\left(  \Phi_{,s}+\Phi\right)  \Phi
_{,s}\Phi}\right)  \partial_{s}+\partial_{\Phi}%
\end{equation}
and
\begin{align}
\bar{X}_{3}  &  =\left(  2\Phi\left(  \frac{1}{\Phi_{,s}}-2\right)
-\frac{2\Phi^{3}}{3\Phi_{,s}^{2}}-\frac{\Phi_{,ss}-3\Phi^{2}-6}{\Phi^{\prime}%
}-\frac{\Phi^{\prime\prime}}{\Phi}\right)  \partial_{s}\nonumber\\
&  -3\left(  \Phi+\Phi^{\prime}\right)  \left(  2+\Phi+\Phi^{\prime}\right)
\partial_{\Phi}.
\end{align}

At this stage, we have now shown that under the various reductions the
symmetries of the Chazy equation change form.\ Moreover, by using the Lie
point symmetry, $X_{2}$, we reduce the Chazy equation to a second-order
ordinary differential equation which admits two nonlocal symmetries. We
increased the order of the equation by one and rewrote it as a third-order
differential equation which admits generalized symmetries. Of course a
solution for \ (\ref{eq.07}) is also a solution of the Chazy Equation,
(\ref{eq.01}).

As we discussed above, the Chazy Equation has been introduced as an equation
which passes the singularity test but with the solution to be given by a left
Laurent expansion and the existence of a natural barrier. In this regard, we
continue our study with the singularity analysis.

\section{Singularity analysis}

For the singularity analysis we apply the ARS algorithm to the Chazy Equation,
(\ref{eq.01}). For the convenience of the reader we present the basic steps of
the ARS algorithm. Firstly we search for the leading-order behaviour,
$y\left(  x\right)  =y_{0}x^{p}$, in (\ref{eq.01}) (without loss of generality
the location of the moveable singularity can be taken as the origin for an
autonomous equation) and we find that $p=-1$ and $y_{0}=-6$, which is also a
special solution of the differential equation. This solution is related to the
existence of the $sl\left(  2,R\right)  $ algebra.

In order to determine the positions of the constants of integration, i.e. the
resonances, we substitute
\begin{equation}
y\left(  x\right)  =-6x^{-1}+mx^{-1+R}%
\end{equation}
into (\ref{eq.01}) and we take the coefficient of $m$. We find the third-order
polynomial
\begin{equation}
R^{3}+6R^{2}+11R+6=0
\end{equation}
the solutions of which provide us with the positions of the constant of
integration. Hence the resonances are%
\begin{equation}
R=-1~,~R=-2~\text{and~}R=-3.
\end{equation}
The negative resonances tell us that the solution is expressed by a Left
Painlev\'{e} Series%
\begin{equation}
y\left(  x\right)  =%
{\displaystyle\sum\limits_{i=1}^{\infty}}
a_{i}x^{-1-i}-6x^{-1} \label{ss.001}%
\end{equation}
in which $a_{1},~a_{2}$ are the two constants of integration. The last step of
the ARS algorithm is to test the consistency of the solution. For that the
expansion (\ref{ss.001}) is substituted in (\ref{eq.01}) and we see that it is
a solution.

However, when the singularity test fails for a differential equation, it does
not mean that the equation is not integrable but that the solution cannot be
described by a Laurent expansion. On the other hand there exists the
possibility that the equation does not admit any singularity, as in the case
of the linear equation%
\begin{equation}
y^{\prime\prime}-y=0 \label{ss.002}%
\end{equation}
which the solution is
\begin{equation}
y\left(  x\right)  =y_{1}e^{x}+y_{2}e^{-x}.
\end{equation}

However, the linear second-order differential equation, (\ref{ss.002}), can be
written as a third-order nonlinear differential equation by the transformation
$x=-\ln\left(  u\left(  v\right)  \right)  ~,~y=\frac{du\left(  v\right)
}{dv}$, that is%
\begin{equation}
u^{2}u_{,v}u_{,vvv}+\left(  u~\left(  u_{v}\right)  ^{2}-u^{2}\right)
u_{,vv}-\left(  u_{,v}\right)  ^{4}=0 \label{ss.003}%
\end{equation}
which can easily be seen by the singularity test. Hence singularity analysis
is coordinate dependent and that also arises from the property that
singularity analysis can be applied only when the equation is written in
rational form.

We now consider the simple coordinate transformation, $y\rightarrow w^{-1}$,
applied to the Chazy Equation, (\ref{eq.01}), and apply the ARS algorithm. We
find that there are two possible leading-order behaviours, namely%
\begin{equation}
w_{1}\left(  x\right)  =w_{0}x^{-1}~\mbox{\rm and}~w_{2}\left(  x\right)
=\bar{w}_{0}x^{-2}%
\end{equation}
where $w_{0}~\mbox{\rm and}~\bar{w}_{0}$ are arbitrary constants. The
application of the ARS algorithm for the leading behaviour $w_{2}\left(
x\right)  $ tells us that the singularity analysis fails. However, the
algorithm succeeds for $w_{1}\left(  x\right)  $ for which the resonances are
\[
R=-1~,~R=0~,~R=1.
\]

The zero resonance tell us that the coefficient of the leading-order term is
arbitrary and it is a constant of integration. Hence the solution is given by
the Right Painlev\'{e} Series%
\begin{equation}
w\left(  x\right)  =w_{0}x^{-1}+w_{1}+w_{2}x+%
{\displaystyle\sum\limits_{i=3}^{\infty}}
w_{i}x^{-1+i}. \label{ss.004}%
\end{equation}
which is a solution for the Chazy Equation without the existence of a natural
barrier. It is obvious that this is a totally different solution from\ that of
(\ref{ss.001}). Last but not least we mention that the consistency test
provides that (\ref{ss.004}) is indeed a solution.

That an equation can possess more than one expansion is easily illustrated by
the Painlev\'{e}-Ince Equation,
\begin{equation}
y^{\prime\prime}+3yy^{\prime3}=0. \label{e1}%
\end{equation}
The singularity analysis of this equation was given in \cite{Lemmer93a}.

The leading-order behaviour is $y(x)=a(x-x_{0})^{-1}$. The nontrivial values
which $a$ can take are $a=1$ and $a=2$. For the former, the resonances are the
generic $-1$ and $2$. For the latter the resonances are the generic $-1$ and
$-2$. Thus the first expansion is a Right Painlev\'{e} Series and the second a
Left Painlev\'{e} Series \cite{Feix97a}.

We continue with the application of the ARS algorithm to the third-order
differential equation (\ref{eq.07}). Recall that this equation is related with
the Chazy Equation under some (nonpoint) generalized transformations.

The application of the ARS algortithm for equation (\ref{eq.07}) \ gives us
the leading-order behaviour with resonances%
\[
\Phi\left(  x\right)  =-6x^{-1},~R=-1~,R=-2~,R=-3,
\]
of which the leading-order behaviour is a special solution of the equation.
Hence the solution is expressed in terms of the Left Painlev\'{e} Series%
\begin{equation}
\Phi\left(  x\right)  =%
{\displaystyle\sum\limits_{i=1}^{\infty}}
\Phi_{i}x^{-1-i}-6x^{-1}.
\end{equation}
However, what it is more interesting is that the leading-order behaviour and
the resonances are exactly the same as that of the Chazy equation.

When we apply the coordinate transformation $\Phi\rightarrow\Psi^{-1}$, the
ARS algorithm gives us the second possible solution%
\begin{equation}
\Psi\left(  x\right)  =\Psi_{0}x^{-1}+\Psi_{1}+\Psi_{2}x+%
{\displaystyle\sum\limits_{i=3}^{\infty}}
\Psi_{i}x^{-1+i}%
\end{equation}
in which the corresponding leading-order term and the resonances are%
\begin{equation}
\Psi\left(  x\right)  =\Psi_{0}x^{-1}~,~R=-1,~R=0~,R=1
\end{equation}
which is again in comparison with solution (\ref{ss.004}). However, what
differs in both cases are the coefficients of the Painlev\'{e} Series. They
are not identical which means that the solution functions are different.
However, they are related with the generalized transformation which transforms
(\ref{eq.07}) \ to the Chazy Equation.

\section{Conclusions}

In this work we studied the Chazy Equation by using two methods for the study
of the integrability of differential equations, the method of group invariant
transformations and that of singularity analysis. These two methods provide
different results. The Chazy Equation, (\ref{eq.01}), admits the elements of
the $sl\left(  2,R\right)  $ algebra as Lie point symmetries. We use the Lie
invariants in order to reduce the differential equation and at the same time
we show how the symmetries changed type, that is, from point symmetries they
became nonlocal symmetries for the reduced equation. However, when we reduced
with the symmetry vector $X_{2}$ and then increased the order of the
second-order reduced equation in order to have an autonomous third-order
ordinary differential equation, we can see that the original symmetries,
$X_{1}$ and $X_{3}$, are now generalized symmetries for the new equation. This
means that there exists a generalized (nonpoint) transformation which connects
that equation to the Chazy Equation.

We continued with the application of the ARS algorithm to the Chazy Equation
and we recovered the well-known solution with negative resonances. However, we
found that the Chazy Equation admits another solution which is given in terms
of a Right Painlev\'{e} Series. However, the most interesting result is that
the new third-order equation which we constructed has also two solutions given
by the singularity analysis in which the leading-order behaviors and the
corresponding resonances are in comparison with those of the Chazy Equation.

The common practice in analysing differential equations using symmetry is to
confine one's attention to point symmetries. Procedurally this is the simpler
approach, but it should always be bourne in mind that point symmetries are a
subset of the whole gamut of symmetries which an equation can possess. Indeed
there are times when it is necessary to move out of the comfort of point
symmetries to reveal the truth. An instance of this can be seen in
\cite{Leach02a}. There nonlocal symmetries were necessary to show the
connection of integrating factors with the presence of symmetry.

In the practice of singularity analysis there has been a long-standing lack of
understanding of the meaning of negative resonances (apart from the generic
$-1$) and yet the explanation is so simple that it is contained in a first
course on complex analysis. The message would seemingly be to think carefully
before condemning an equation to the realm of nonintegrability.

Finally, it will be of special interest to follow the approach of \cite{abl1}
and study the relation for the solution of equation\ (\ref{eq.07}) with that
of the Schwarz triangle functions, and also study the effects of nonlocal
transformation in the solution of the original equation.

\begin{acknowledgments}
AP acknowledges the financial support of FONDECYT grant no. 3160121 and thanks
the Durban University of Technology and the University of the Witwatersrand
for the hospitality provided while this work was performed. SJ would like to
acknowledge the financial support from the National Research Foundation of
South Africa with grant number 99279. PGLL acknowledges the National Research
Foundation of South Africa and the University of KwaZulu-Natal for financial
support. The authors acknowledges support from the DST-NRF Centre of
Excellence in Mathematical and Statistical Sciences (CoE-MaSS). Special thanks
to Caryn McNamara and Shashilan Singh of COE-Mass. Opinions expressed and
conclusions arrived are those of the authors and are not to be construed as
being those of any of the institutions above.
\end{acknowledgments}

\end{document}